\documentclass[12pt]{amsart}% other options are: book, amsart,
\usepackage{amsmath,amsthm,amsfonts,amssymb}% standard options;
\usepackage[% headsep=1cm,
margin=2.5cm]{geometry}
\usepackage[usenames,dvipsnames]{color}% allows manual adjustments to
\usepackage{mathrsfs}% enables mathscr (and other font stuff...?)
\usepackage[colorlinks,allcolors=blue]{hyperref}% adds hyperlinks to
\usepackage{enumitem}% for controlling indents in enumerated
\setlist[enumerate]{%
  leftmargin=*,%
  label=(\roman*)} % apparently these have to go after the label
\setlist[itemize]{leftmargin=*}
\usepackage[mathscr]{eucal}% Greg's suggestion for the trees
\usepackage{multirow, array}% for the new table

\usepackage{tikz,xcolor}
\usetikzlibrary{positioning, backgrounds}

\usepackage{fancyvrb} %%--FANCY-VERBATIM
\RecustomVerbatimEnvironment{Verbatim}{Verbatim}%
{xleftmargin=0pt,baselinestretch=0.85}

\usepackage{microtype}
\usepackage[yyyymmdd]{datetime}% international standard
\usepackage{fancyhdr}% the standard package for formatting headers and
\newtheorem{theorem}{Theorem}[section]
\newtheorem{lemma}[theorem]{Lemma}
\newtheorem{proposition}[theorem]{Proposition}

\theoremstyle{definition}%
\newtheorem{example}[theorem]{Example}
\newtheorem{remark}[theorem]{Remark}

\DeclareMathOperator{\hilb}{Hilb}% notation for the Hilbert scheme
\DeclareMathOperator{\PP}{\mathbb{P}}% projective space
\DeclareMathOperator{\bbA}{\mathbb{A}}% affine space
\DeclareMathOperator{\Proj}{Proj}% Proj of a graded ring
\DeclareMathOperator{\Spec}{Spec}% Spec of a ring
\DeclareMathOperator{\rank}{rk}% the regularity

\DeclareMathOperator{\GL}{GL}% 
\DeclareMathOperator{\Hom}{Hom}% 

\newcommand{\hf}{\mathsf{h}}% the Hilbert function
\newcommand{\hp}{\mathsf{p}}% the Hilbert polynomial
\newcommand{\HG}{\operatorname{Hilb}^{4d}(\PP^3)}

\newcommand{\kk}{\mathbb{K}}% the ground field
\newcommand{\NN}{\mathbb{N}}% the nonnegative integers
\newcommand{\QQ}{\mathbb{Q}}% the rational field
\newcommand{\ZZ}{\mathbb{Z}}% the integers
\DeclareMathOperator{\reg}{reg}% the regularity
\DeclareMathOperator{\gin}{gin}% generic initial ideal
\DeclareMathOperator{\ini}{in}% initial ideal
\newcommand{\llrr}[1]{ \langle #1 \rangle } % short for \langle \rangle
\begin{document}

% --------------------------------------------------------------------
% Abstract
% --------------------------------------------------------------------
\begin{abstract}
  We study the component structures of some standard-graded Hilbert
  schemes closely related to a Hilbert scheme of curves studied by
  Gotzmann.  In particular, we encounter examples of singular
  lex-segment points lying on two and three irreducible components.
  We find further singular lex-segment points at nearby Hilbert
  schemes.  We conclude by showing that the analogous example at the
  Hilbert scheme of twisted cubics also has a singular lex-segment
  point.
\end{abstract}
% --------------------------------------------------------------------

% --------------------------------------------------------------------
% Title Page
% --------------------------------------------------------------------
\title{Some Singular Lex-Segments}
\author{Andrew~P.~Staal}%
\email{andrew.staal@uwaterloo.ca}
\address{ Department of Pure Mathematics \\
  University of Waterloo \\
  200 University Avenue West \\
  Waterloo, Ontario, Canada, N2L 3G1 }%
\date{}
\maketitle
% --------------------------------------------------------------------

% --------------------------------------------------------------------
% Section 0 -- Introduction
% --------------------------------------------------------------------
\section{Introduction}

In \cite{Staal--2020}, we provide a classification of Hilbert schemes
$\hilb^{\hp}(\PP^n)$, parametrizing closed subschemes of projective
space $\PP^n$ with a specified Hilbert polynomial $\hp$, that have
unique Borel-fixed points.  These Hilbert schemes are smooth and our
classification reveals the existence of natural probability
distributions where smooth Hilbert schemes occur with probability
greater than $0.5$.  Based on our underlying ``geography'' of Hilbert
schemes, further smooth examples are promptly found, leaving open the
challenge of understanding whether all smooth Hilbert schemes fall
into known classes.  This is achieved by Skjelnes--Smith in
\cite{Skjelnes--Smith--2020}, where our classification comprises two
of the main classes of smooth Hilbert schemes.

Multigraded Hilbert schemes $\hilb^{\hf}(S)$, parametrizing
homogeneous ideals in a multigraded polynomial ring $S$ with a
specified Hilbert function $\hf$, generalize Hilbert schemes by using
a broader notion of degree.  Introduced in
\cite{Haiman--Sturmfels--2004}, these quasi-projective schemes include
Hilbert schemes of points in affine space and toric Hilbert schemes,
among other examples.  Their geometry is not well-understood.  All
pathologies that occur for Hilbert schemes also occur for multigraded
Hilbert schemes, but further complications arise: toric Hilbert
schemes can be disconnected \cite{Santos--2005} and lexicographic
points do not exist for general multigradings.

A natural question is: To what extent does the classification of
smooth Hilbert schemes extend to multigraded Hilbert schemes?  For
example, multigraded Hilbert schemes parametrizing admissible ideals
in two variables are nonsingular \cite{Maclagan--Smith--2010}.  A key
fact used in these classifications is that lex-segment, or lex-most,
ideals are smooth \cite{Reeves--Stillman--1997}.  Thus, an important
first step towards the general case is to understand the geometry of
multigraded Hilbert schemes at lex-segment ideals, when they exist.
In this paper, we study the component structures of some explicit
standard-graded Hilbert schemes closely related to Hilbert schemes.

Let $S := \kk[x_0, x_1, \ldots, x_n]$, where $\kk$ is a field of
characteristic $0$.  We prove the following.

\begin{theorem}
  Let $M^{\hf} = \bigcup_{i=0}^s M^{k_i}_i$ be the irreducible
  decomposition of the standard-graded Hilbert scheme $M^{\hf} :=
  \hilb^{\hf}(S)$, where $M^{k_i}_i$ has dimension $k_i$.  Let $I_0$
  be the lex-segment ideal with Hilbert function $\hf$.  In the
  following cases, let $\hf_i$ be the Hilbert function of the
  saturated Borel-fixed ideal having Hilbert polynomial $\hp(d)$ and
  regularity $\reg(I_0) - i$.  We have
  \begin{itemize}
  \item for $n=3$ and Hilbert polynomial $\hp(d) = 4d$:
    \begin{itemize}
    \item $M^{\hf_1} = M^{24}_0 \cup M^{22}_1$ with $I_0 \in M^{24}_0
      \cap M^{22}_1$,
    \item $M^{\hf_2} = M^{28}_0 \cup M^{25}_1 \cup M^{23}_2 \cup
      M^{16}_3$ with $I_0 \in M^{28}_0 \cap M^{25}_1 \cap M^{23}_2$;
    \end{itemize}

  \item for $n=3$ and Hilbert polynomial $\hp(d) = 4d-1$:
    \begin{itemize}
    \item $M^{\hf_1} = M^{22}_0 \cup M^{20}_1$ with $I_0 \in M^{22}_0
      \cap M^{20}_1$;
    \end{itemize}

  \item for $n=2$ and Hilbert polynomial $\hp(d) = 4$:
    \begin{itemize}
    \item $M^{\hf_1} = M^{8}_0 \cup M^{8}_1$ with $I_0 \in M^{8}_0
      \cap M^{8}_1$;
    \end{itemize}

  \item for $n=3$ and Hilbert polynomial $\hp(d) = 3d+1$:
    \begin{itemize}
    \item $M^{\hf_1} = M^{17}_0 \cup M^{15}_1 \cup M^{12}_2$ with $I_0
      \in M^{17}_0 \cap M^{15}_1$.
    \end{itemize}
  \end{itemize}
\end{theorem}

In stark constrast to $\hilb^{\hp}(\PP^n)$, these lex-segment points
all lie in nonempty intersections of irreducible components.  The
examples we study are roughly organized, first, by admissible Hilbert
polynomial according to the Macaulay tree \cite[\S 2]{Staal--2020}
and, second, according to Hilbert functions that appear naturally when
studying (classical) Hilbert schemes; cf.\ \cite{Gotzmann--2008}.
Irreducible components then emerge by considering ideals of fixed
regularity.

A first example of a singular lex-segment ideal is given in
\cite{Ramkumar--Sammartano--2020}.  Our examples broaden the class of
known singular lex-segments and raise many interesting questions: Are
all lex-segment points of this kind singular?  Do they always lie in
the largest numbers of irreducible components?  Are these numbers of
irreducible components containing lex-segment points unbounded when we
vary the Hilbert schemes in natural ways?  Does Murphy's Law hold at
lex-segment points?  Does the classification of smooth Hilbert schemes
extend ``up to'' these examples with singular lex-segments?

In subsequent sections, we provide a more detailed analysis of the
component structures of the standard-graded Hilbert schemes listed
above, as well as pointing out some nonsingular standard-graded
Hilbert schemes along the way.

\subsection*{Conventions.}

Throughout, $\kk$ is a field of characteristic $0$, $\NN$ denotes the
nonnegative integers, and (after \S \ref{BF}) $S := \kk[x_0, x_1,
  \ldots, x_n]$ is a standard $\ZZ$-graded polynomial ring over $\kk$,
endowed with the reverse lexicographic monomial ordering unless
otherwise stated.

\subsection*{Acknowledgments}

We thank Chris Brav for helpful suggestions and encouragement, and
Matt Satriano for suggestions for improving the presentation.  This
research was supported by a Geometry and Topology postdoctoral
fellowship at the University of Waterloo.
% --------------------------------------------------------------------

% --------------------------------------------------------------------
% Section 1 -- Useful Tools
% --------------------------------------------------------------------
\section{Background and Useful Tools}

In this section, we highlight some basic facts and useful tools for
studying Hilbert schemes.

\subsection{Basic Facts}
\label{BF}

Let $R$ be a unital commutative ring and $S := R[x_0, x_1, \ldots,
  x_n]$ be the polynomial ring over $R$.  A \emph{\bfseries
multigrading} on $S$ is a monoid homomorphism $\deg \colon \NN^{n+1}
\to A$, where $A$ is an abelian group; the submonoid $A_+ :=
\deg(\NN^{n+1})$ is assumed to generate $A$ as a group.  A
multigrading determines a decomposition $S = \bigoplus_{d \in A} S_d$
into multigraded components $S_d := \bigoplus_{\deg(u) = d} R x^u$,
where $x^u := x_0^{u_0} x_1^{u_1} \cdots x_n^{u_n}$.  We use the
standard grading $\deg \colon \NN^{n+1} \to \ZZ$, $(u_0, u_1, \ldots,
u_n) \mapsto u_0 + u_1 + \cdots + u_n$, but refer to some facts that
apply more generally.  A multigrading is \emph{\bfseries positive} if
$S_0 = R$---the standard grading is positive.  An ideal $I \subseteq
S$ is \emph{\bfseries admissible} if it is homogeneous with respect to
the multigrading and $S_d/I_d$ is a locally free $R$-module of finite
rank, for all $d \in A$.  When $S$ is positively multigraded, this is
equivalent to $S/I$ being flat over $R$ \cite[\S
  18.5]{Miller--Sturmfels--2005}.  For an admissible ideal $I$, the
\emph{\bfseries Hilbert function} $\hf_{S/I} \colon A \to \ZZ$ of
$S/I$ is defined by $\hf_{S/I}(d) := \rank_{\kk} (S_d/I_d)$, for all
$d \in A$.  It is common to refer to $\hf_{S/I}$ as the Hilbert
function of $I$ and use the notation $\hf_I$, however, working with
the standard grading, we refer to $\hf_I(d) := \hf_S(d) -
\hf_{S/I}(d)$ as the Hilbert function of $I$.  Every such
standard-graded $I$ has a \emph{\bfseries Hilbert polynomial}
$\hp_{S/I}$, that is, a polynomial $\hp_{S/I}$ with $\QQ$-coefficients
such that $\hf_{S/I}(d) = \hp_{S/I}(d)$ for $d \gg 0$; see
\cite[Theorem~4.1.3]{Bruns--Herzog--1993}.

Admissible ideals are parametrized by multigraded Hilbert schemes,
whose existence and fundamental properties are expounded by
Haiman--Sturmfels in \cite{Haiman--Sturmfels--2004}.  We collect some
key facts about multigraded Hilbert schemes in the following theorem.

\begin{theorem}
  \label{HS}
  Let $S$ be multigraded and $\hf \colon A \to \NN$ be a function
  supported on $A_+$.
  \begin{enumerate}
    \item The admissible ideals in $S$ with Hilbert function
      $\hf_{S/I} = \hf$ are parametrized by a quasi-projective scheme
      $\hilb^{\hf}(S)$ called the {\bfseries multigraded Hilbert
      scheme}.
    \item If the multigrading on $S$ is positive, then
      $\hilb^{\hf}(S)$ is projective.
    \item If $R = \kk$ is a field, then the Zariski tangent space to
      $\hilb^{\hf}(S)$ at a $\kk$-point $I \subseteq S$ equals
      $\Hom_S(I, S/I)_0$, the degree-$0$ part of the multigraded
      $S$-module $\Hom_S(I, S/I)$.
  \end{enumerate}
\end{theorem}

\begin{proof}
  See Section 1 of \cite{Haiman--Sturmfels--2004}.
\end{proof}

\begin{example}  Many common parameter spaces are multigraded Hilbert
  schemes.
  \begin{enumerate}
  \item The Hilbert scheme $\hilb^{\hp}(\PP^n)$ is isomorphic to the
    multigraded Hilbert scheme $\hilb^{\hf}(S)$, where $\hf(d) =
    \hp(d)$, for $d \ge d_0$, and $\hf(d) = \binom{n+d}{n}$, for $d <
    d_0$ (meaning $0$ when $d < 0$); here $d_0$ is the Gotzmann number
    of $\hp$ \cite[Lemma 4.1]{Haiman--Sturmfels--2004}.
  \item If $A = 0$ is trivial, then $\hilb^{\hf}(S)$ is the Hilbert
    scheme of $\hf(0)$ points in $\bbA^{n+1}$.
  \item Let the multigrading be positive with $A \cong \ZZ^r$ and set
    $\hf(d) = 1$, for $d \in A_+$, and $\hf(d) = 0$ otherwise.  Then
    $\hilb^{\hf}(S)$ is the toric Hilbert scheme of Peeva--Stillman
    \cite{Peeva--Stillman--2002}.
  \end{enumerate}
\end{example}

When $S$ is standard-graded and $I$ is admissible, we refer to
$\hilb^{\hf_I}(S)$ to mean $\hilb^{\hf_{S/I}}(S)$.  We also use the
notation $M^{\hf} := \hilb^{\hf}(S)$, when $S$ is clear from context.

The saturation of an ideal $I \subseteq S$ is the ideal
\[
I^{sat} := \left( I \colon \mathfrak{m}^{\infty} \right) = \bigcup_{k
  \ge 1} \left\{ f \in S \mid f \mathfrak{m}^k \subseteq I \right\},
\]
where $\mathfrak{m} := \llrr{x_0, x_1, \ldots, x_n}$ is the irrelevant
ideal of $S$.  The inclusion $I \subseteq I^{sat}$ always holds; if
$S$ is standard-graded and $I$ is homogeneous, then $I^{sat}$ is
homogeneous, $I_d \subseteq I^{sat}_d$ for all $d \in \ZZ$, and $I_d =
I^{sat}_d$ for $d \gg 0$.

\subsection{Initial Ideals}

Let $S := \kk[x_0, x_1, \ldots, x_n]$ be standard-graded and let $>$
be a monomial ordering with $x_0 > x_1 > \cdots > x_n$; we assume a
familiarity with monomial orderings and initial ideals, as provided in
\cite{Cox--Little--O'Shea--2015}.  The following result is often used
in the study of the component structure of the Hilbert scheme.

\begin{theorem}
  Given an ideal $I \subseteq S$ and a monomial ordering $>$ on $S$,
  there exists a family $J = \{ J_a \} \subseteq S[a]$ of ideals such
  that
  \begin{enumerate}
  \item $J_1 = I$ and $J_a \cong I$ by a scaling of variables, for $a
    \neq 0$, and
  \item $\lim_{a \to 0} J_a = \ini I$.
  \end{enumerate}
\end{theorem}

\begin{proof}
  See \cite[\S 2.12]{Bayer--1982}, \cite[Theorem 3]{Reeves--1992}, or
  \cite[\S 15.8]{Eisenbud--1995}.
\end{proof}

Here $\ini I = \ini_> (I) \subseteq S$ is the initial ideal of $I$
with respect to $>$.  To study standard-graded Hilbert schemes, we
assume $I$ is homogeneous and work directly with the family $J$ rather
than with the flat family $\Proj S[a]/J \to \Spec \kk[a] \cong
\bbA^1_{\kk}$.  The theorem says that given any homogeneous ideal $I$,
there is a morphism $\Spec \kk[a] \to \hilb^{\hf_I}(S)$ with fibre $I$
at $a=1$ and fibre $\ini_> (I)$ at $a=0$.  Consider the graded
$S[a]$-module $\Hom_{S[a]}(J, S[a]/J)$, where $J \subseteq S[a]$ is
homogeneous (and $\deg a = 0$).  Then $\Hom_{S[a]}(J, S[a]/J)_0$ is a
$\kk[a]$-module and another useful result is the following.

\begin{lemma}
  The dimension of $\Hom_S(J_a, S/J_a)_0$ is upper semi-continuous on
  $\Spec \kk[a]$.
\end{lemma}

\begin{proof}
  By Theorem \ref{HS}(iii), the $\kk$-vector space $\Hom_S(J_a,
  S/J_a)_0$ is the Zariski tangent space $T_{J_a} M$ to the
  standard-graded Hilbert scheme $M := \hilb^{\hf_{J_0}}(S)$ at the
  point $J_a$ and, more generally, the dimension of $T_{I} M$ is upper
  semi-continuous for $I$ in $M$ \cite[\S 6.6]{Goertz--Wedhorn--2010}.
\end{proof}

In practice, these two results tell us that the initial ideal lies on
the same component of the standard-graded Hilbert scheme as the ideal
itself and that the dimension of the tangent space cannot decrease
when specializing to the initial ideal.

The action of an element $\gamma \in \GL_{n+1}(\kk)$ on $S$ is defined
by $\gamma \cdot x_j := \sum_{i=0}^n \gamma_{ij} x_i$.  The
Borel-subgroup of $\GL_{n+1}(\kk)$ consists of the upper-triangular
matrices.  The \emph{\bfseries generic initial ideal} $\gin I = \gin_>
(I) \subseteq S$ of a homogeneous ideal $I \subseteq S$ is the initial
ideal of a generic change-of-coordinates of $I$.  Specifically, we
have the following.

\begin{theorem}[Galligo]
  If $I \subseteq S$ is a homogeneous ideal, then there exists a
  Zariski open subset $U \subseteq \operatorname{GL}_{n+1}(\kk)$ such
  that $\ini \gamma I = \ini \gamma' I$, for all $\gamma, \gamma' \in
  U$.  Moreover, this common initial ideal, called $\gin I$, is
  Borel-fixed.
\end{theorem}

A monomial ideal $I \subseteq S$ is \emph{\bfseries strongly stable}
if, for all monomials $m \in I$, for all $x_j$ dividing $m$, and for
all $i < j$, we have $x_j^{-1} m x_i \in I$.  In characteristic $0$,
this condition is equivalent to being \emph{\bfseries Borel-fixed}.
(In positive characteristic, the condition is more involved.)  The
\emph{\bfseries Borel partial ordering} on the monomials of $S$ is the
transitive closure of the following: for each monomial $m$, for all
$x_j$ dividing $m$, and all $i < j$, the monomial $x_j^{-1} m x_i$ is
larger than $m$.

\begin{proof}
  The original is \cite{Galligo--1974}, but we are interested in the
  generalization by Bayer--Stillman \cite[Proposition
    1]{Bayer--Stillman--1987b}, \cite[\S 15.9]{Eisenbud--1995}.
\end{proof}

Combining Galligo's theorem with the degeneration to the initial ideal
shows that every irreducible component and intersection of irreducible
components of a standard-graded Hilbert scheme $\hilb^{\hf}(S)$
contains a Borel-fixed point.

\begin{example}[Lex-segment ideals]
  The most important Borel-fixed ideals are the lex-segment ideals.
  These are defined with respect to the \emph{\bfseries lexicographic
  ordering}, in which $x^u > x^v$ if $u_i > v_i$, where $i$ is minimal
  such that $u_i \neq v_i$.  Given a homogeneous ideal $I \subseteq S$
  with Hilbert function $\hf = \hf_I$, the \emph{\bfseries lex-segment
  ideal} with Hilbert function $\hf$ is the ideal $L$ whose degree-$d$
  component $L_d$ is spanned by the $\hf(d)$ greatest monomials in
  lexicographic order; that $L$ is an ideal is a theorem of Macaulay
  \cite{Macaulay--1927}, \cite[\S 2.4]{Miller--Sturmfels--2005}.  The
  corresponding nonsingular point $\Proj S/L \in \hilb^{\hp}(\PP^n)$
  is called the lexicographic point \cite{Reeves--Stillman--1997}.
\end{example}

Another useful monomial ordering is the \emph{\bfseries reverse
lexicographic ordering}, in which $x^u > x^v$ if $u_j < v_j$, where
$j$ is maximal such that $u_j \neq v_j$.  Despite the importance of
lexicographic order, the reverse lexicographic ordering is optimized
to computing Gr\"{o}bner bases \cite{Bayer--Stillman--1987b}.  We use
the reverse lexicographic ordering unless specified otherwise.

\subsection{Deformation Theory}

Although nearly superfluous for gaining a basic understanding of the
component structures described in the following sections, deformation
theory offers a powerful experimental tool for probing the
(analytically) local geometry of Hilbert schemes.  This is done via
the \emph{power series ansatz}, which is implemented
\cite{Ilten--2012} in \emph{Macaulay2} \cite{Grayson--Stillman}.  The
power series ansatz is the procedure that lifts an ideal,
degree-by-degree, over a ring of deformation parameters (of dimension
equal to that of the Zariski tangent space at the ideal in the Hilbert
scheme or a local analogue).  At each stage, corrections might need to
be made to eliminate the obstructions to lifting the ideal.  The
implementation is described in \cite{Ilten--2012}, an introduction to
the process is given in \cite{Stevens--1995}, and a classic example
done by hand is provided in \cite[\S 5]{Piene--Schlessinger--1985} (a
detailed exposition in \cite{Staal--2014} summarizes the basics and
follows concrete examples).  The ansatz may not terminate, but when it
does it provides a local view of the multigraded Hilbert scheme.  As
in the case of twisted cubics, it can also provide a local view of the
Hilbert scheme, when the Comparison Theorem holds \cite[\S
  3]{Piene--Schlessinger--1985}.

\begin{example}
  \label{defeg}
  We demonstrate how to quickly produce a description of the local
  geometry of the scheme $\hilb^{\hf_1}(S)$ at the lex-segment ideal
  $I_0' := \llrr{x^2, xy, xz, xt^2, y^4, y^3z}$ in $S := \kk[x,y,z,t]$
  with $\hf_1 = \hf_{I_0'} = (0, 0, 3, 10, 22, 40, 65, \ldots)$; the
  component structure of this scheme is described in
  Theorem~\ref{twist}.  Note that $I_0'$ has Hilbert polynomial
  $\hp_{S/I_0'}(d) = 3d+1$.
\begin{Verbatim}[fontsize=\small]
Macaulay2, version 1.17.2.1
with packages: ConwayPolynomials, Elimination, IntegralClosure, InverseSystems,
               LLLBases, MinimalPrimes, PrimaryDecomposition, ReesAlgebra,
               Saturation, TangentCone
\end{Verbatim}

\begin{Verbatim}[fontsize=\small]
i1 : S = QQ[x,y,z,t];
\end{Verbatim}
\begin{Verbatim}[fontsize=\small]
i2 : I0' = ideal(x^2, x*y, x*z, x*t^2, y^4, y^3*z);
\end{Verbatim}
\begin{Verbatim}[fontsize=\small]
o2 : Ideal of S
\end{Verbatim}
The next line computes the dimension of the tangent space $T_{I_0'}
\hilb^{\hf_1}(S)$.
\begin{Verbatim}[fontsize=\small]
i3 : n = rank source ambient basis(0, Hom(I0', S/I0'))
\end{Verbatim}
\begin{Verbatim}[fontsize=\small]
o3 = 18
\end{Verbatim}
As explained in \S 5, we expect $I_0'$ to lie in two irreducible
components, of dimensions $17$ and $15$.  To see this computationally,
we proceed as follows.
\begin{Verbatim}[fontsize=\small]
i4 : needsPackage "VersalDeformations"
--loading configuration for package "VersalDeformations" from file /Users/ ...
\end{Verbatim}
\begin{Verbatim}[fontsize=\small]
o4 = VersalDeformations
\end{Verbatim}
\begin{Verbatim}[fontsize=\small]
o4 : Package
\end{Verbatim}
\begin{Verbatim}[fontsize=\small]
i5 : time (F,R,G,C) = localHilbertScheme(gens I0', Verbose => 2); use S;
Calculating first order deformations and obstruction space
Calculating first order relations
Calculating standard expressions for obstructions
Starting lifting
Order 2
Order 3
Order 4
Order 5
Order 6
Order 7
Order 8
Order 9
Order 10
Order 11
Order 12
Order 13
Solution is polynomial
     -- used 1.02269 seconds
\end{Verbatim}
Now we extract the obstruction ideal.
\begin{Verbatim}[fontsize=\small]
i7 : ob = ideal sum G
\end{Verbatim}
\begin{Verbatim}[fontsize=\small]
o7 = ideal (- u  u   - u  u  , u  u   - u  u  , u u   - 2u  u  )
               10 18    16 18   11 18    17 18   9 18     14 18
\end{Verbatim}
\begin{Verbatim}[fontsize=\small]
o7 : Ideal of S[u ..u  ]
                 1   18
\end{Verbatim}
We see by inspection that the ideal is as expected from
Theorem~\ref{twist}, with primary decomposition $\llrr{ u_{10}+u_{16},
  u_{11}-u_{17}, u_9-2u_{14}} \cap \llrr{u_{18}} \subset \QQ[u_1, u_2,
  \ldots, u_{18}]$.  The \emph{Macaulay2} commands to see this are as
follows, first removing the variables $x,y,z,t$ from the ambient
obstruction ring.
\begin{Verbatim}[fontsize=\small]
i8 : Ob = QQ[u_1..u_n]; phi = map(Ob, ring ob, vars Ob); ob = phi(ob)
\end{Verbatim}
\begin{Verbatim}[fontsize=\small]
o9 : RingMap Ob <--- S[u ..u  ]
                        1   18
\end{Verbatim}
\begin{Verbatim}[fontsize=\small]
o10 = ideal (- u  u   - u  u  , u  u   - u  u  , u u   - 2u  u  )
                10 18    16 18   11 18    17 18   9 18     14 18
\end{Verbatim}
\begin{Verbatim}[fontsize=\small]
o10 : Ideal of Ob
\end{Verbatim}
\begin{Verbatim}[fontsize=\small]
i11 : time pd = primaryDecomposition ob; #pd
     -- used 0.0389783 seconds
\end{Verbatim}
\begin{Verbatim}[fontsize=\small]
o12 = 2
\end{Verbatim}
\begin{Verbatim}[fontsize=\small]
i13 : for i in pd do << dim i << endl;
17
15
\end{Verbatim}
\begin{Verbatim}[fontsize=\small]
i14 : exit
\end{Verbatim}
\begin{Verbatim}[fontsize=\small]
Process M2 finished
\end{Verbatim}
Hence, we see irreducible components of dimensions $17$ and $15$
meeting at the lex-segment ideal $I_0'$.  Similar analyses work at
other points and on our examples in the coming sections.
\end{example}
% --------------------------------------------------------------------

% --------------------------------------------------------------------
% Section 2 -- The scheme Hilb^{4t}(\PP^3)
% --------------------------------------------------------------------
\section{Standard-graded Hilbert Schemes and an Example of Gotzmann}

In \cite{Gotzmann--2008}, Gotzmann studies the Hilbert scheme $H :=
\HG$ via its regularity stratification, showing that $H = H_{AV} \cup
H_{RS}$ consists of two irreducible components.  The component
$H_{AV}$ is $16$-dimensional and has general point defined by an
arithmetically Cohen--Macaulay curve in $\PP^3$.  It is generically
nonsingular \cite{Ellingsrud--1975}.  The lexicographic component
$H_{RS}$ is $23$-dimensional, has general member defined by a
degenerate quartic curve union two points, and is also generically
nonsingular \cite{Reeves--Stillman--1997}.  The components intersect
transversally and are rational; this is easily seen using the power
series ansatz \cite[\S 1]{Stevens--1995}, implemented in
\emph{Macaulay2} \cite{Ilten--2012}.

We exhibit some interesting properties of standard-graded Hilbert
schemes closely related to $H$.  In particular, these standard-graded
Hilbert schemes provide new examples of singular lex-segment ideals,
lying in multiple irreducible components.

\subsection{Setup and a First Example}

To begin, we list the possible Hilbert functions of saturated ideals
defining points of $H$.  In fact, it suffices to list the Hilbert
functions of saturated Borel-fixed ideals with Hilbert polynomial
$\hp(d) = 4d$.  Letting $S := \kk[x,y,z,t]$ denote the homogeneous
coordinate ring of $\PP^3$, these are as follows:
\begin{align*}
  I_0 &= \llrr{ x,y^5,y^4z^2 }, &\hf_0 &= (0, 1, 4, 10, 20, 36, 60,
  \ldots), \\
  I_1 &= \llrr{ x^2, xy, xz, y^5, y^4z }, &\hf_1 &= (0, 0, 3, 9, 19, 36,
  60, \ldots), \\
  I_2 &= \llrr{ x^2, xy, xz^2, y^4 }, &\hf_2 &= (0, 0, 2, 8, 19, 36, 60,
  \ldots), \\
  I_3 &= \llrr{ x^2, xy, y^3 }, &\hf_3 &= (0, 0, 2, 8, 19, 36, 60, \ldots),
\end{align*}
where $\hf_i := \hf_{I_i}$ and $I_i$ has regularity $6-i$.  The list
of ideals can be derived using Reeves' algorithm \cite{Reeves--1992,
  Moore--Nagel--2014}.  Suppose $X \subset \PP^3$ is a closed
subscheme with Hilbert polynomial $4d$ and let $I_X \subset S$ denote
the saturated ideal of $X$, i.e.\ $I_X := \bigoplus_{d \ge 0}
H^0(\PP^3, \mathscr I_X(d))$ where $\mathscr I_X$ is the ideal sheaf
of $X$.  By Galligo's theorem, $\gin I_X$ is Borel-fixed.  Further,
under the reverse lexicographic ordering, $\gin I_X$ has regularity
equal to that of $I_X$ \cite[\S 2]{Bayer--Stillman--1987} and is
saturated \cite[Theorem~2.30]{Green--2010}, so is in the given list.
In particular, the Hilbert function of $I_X$ equals that of $\gin I_X$
and is listed.

The first standard-graded Hilbert scheme related to $H$ that we
examine parametrizes all homogeneous ideals $J \subset S$ with Hilbert
function $\hf_0 = (0, 1, 4, 10, 20, 36, 60, \ldots)$, i.e.\ equal to
that of the lexicographic ideal $I_0$.  We denote this scheme by
$M^{\hf_0}$ and prove the following.

\begin{proposition}
  \label{h0}
  The scheme $M^{\hf_0}$ is nonsingular and irreducible of dimension
  $21$.
\end{proposition}

\begin{proof}
  Let $J \subset S$ be a point of $M^{\hf_0}$.  Inspecting the Hilbert
  function, one finds $J = \llrr{\ell, f, g}$ has generators $\ell \in
  S_1$, $f = \ell' h \in S_5$, and $g = q h \in S_6$, where $\ell' \in
  S_1 \setminus \kk \ell$, $h \in S_4 \setminus S_3 \ell$, and $q \in
  S_2 \setminus S_1 \llrr{\ell, \ell'}$ \cite[\S 2.4]{Gotzmann--2008}.
  Counting parameters shows this is a $21$-dimensional family, of
  which $I_0$ is a member.  On the other hand, $\gin J$ is Borel-fixed
  with Hilbert function $\hf_0$, which necessitates $\gin J = I_0$
  (see Lemma~\ref{nonsat} below).  The lex-segment point $I_0$ is
  nonsingular, as one verifies as in Example \ref{defeg} that
  $\dim_{\kk} T_{I_0} M^{\hf_0} = \dim_{\kk} \Hom_S(I_0, S/I_0)_0 =
  21$.  By upper semi-continuity of $\dim_{\kk} T_I M^{\hf_0}$, the
  point $J$ is also nonsingular.
\end{proof}

This simple proof generalizes to standard-graded Hilbert schemes with
unique lex-segment ideals that are nonsingular.  That is, the
following holds, for any $S := \kk[x_0, x_1, \ldots, x_n]$.

\begin{proposition}
  \label{h0'}
  Let $M := \hilb^{\hf}(S)$ be a standard-graded Hilbert scheme.  If
  $M$ has a unique Borel-fixed point $I \subseteq S$ and $I$ is
  nonsingular, then $M$ is nonsingular and irreducible.
\end{proposition}

\begin{proof}
  For every $J \in M$, we have $\gin J = I$ and so upper
  semi-continuity of the dimension of the tangent space proves the
  claim.
\end{proof}

\subsection{A Lex-Segment Ideal on Two Components}

Next we study the standard-graded Hilbert scheme $M^{\hf_1}$
parametrizing ideals of $S := \kk[x,y,z,t]$ with Hilbert function
$\hf_1$.  First, we compile a list of Borel-fixed ideals with Hilbert
function $\hf_1$.  To do so, we use the following straight-forward
fact.

\begin{lemma}
  \label{nonsat}
  Suppose that $I \subset S$ is a non-saturated Borel-fixed ideal and
  $d$ is a positive integer.  Then $I^{sat}$ is Borel-fixed and $I_d$
  is obtained from $I^{sat}_d$ by successively removing monomials that
  are minimal in Borel partial order.
\end{lemma}

\begin{proof}
  Let $f \in I^{sat}$, so that all $m \in \mathfrak{m}^k$ satisfy $fm
  \in I$, for some $k > 0$.  If $x_j \vert f$ and $i < j$, then
  $(x_j^{-1}fx_i)m = x_j^{-1}(fm)x_i \in I$ implies $x_j^{-1}fx_i \in
  I^{sat}$, i.e.\ $I^{sat}$ is Borel-fixed.  The containment $I_d
  \subseteq I^{sat}_d$ holds for all $I$ and $d$.  Suppose that
  $\hf_I(d) < \hf_{I^{sat}}(d)$, so that a monomial basis of $I_d$ is
  obtained from a monomial basis of $I^{sat}_d$ by removing monomials.
  No removed monomial is greater than a monomial of $I_d$, in the
  Borel partial order.  As the set of removed monomials is finite, its
  elements can be enumerated by selecting a minimal one and repeating
  on those that remain.
\end{proof}

Let $I$ be Borel-fixed with Hilbert function $\hf_1 = (0, 0, 3, 9, 19,
36, 60, \ldots)$.  Then either $I^{sat} = I_0$ or $I^{sat} = I_1$, as
$\hf_I(2) = 3 > 2 = \hf_2(2) = \hf_3(2)$.  Having $I^{sat} = I_1$
would imply $I = I_1$, so suppose that $I^{sat} = I_0$.  Then $I$ is
obtained from $I_0$ by removing a single monomial in each degree from
$1$ to $4$.  There is a unique way to do this: remove $x$, $xt$,
$xt^2$, and $xt^3$.  In other words, we obtain the list
\begin{align*}
  I_0' &= \llrr{ x^2, xy, xz, xt^4, y^5, y^4z^2 }, \\
  I_1 &= \llrr{ x^2, xy, xz, y^5, y^4z },
\end{align*}
of Borel-fixed ideals with Hilbert function $\hf_1$.  The ideal $I_0'$
is the lex-segment ideal for $\hf_1$.

Given this complete list of Borel-fixed ideals on $M^{\hf_1}$, we
study the component structure of the standard-graded Hilbert scheme
$M^{\hf_1}$.  In particular, we discover a standard-graded lex-segment
ideal contained in the intersection of the two irreducible components
of $M^{\hf_1}$.

\begin{theorem}
  \label{2comps}
  The scheme $M^{\hf_1}$ consists of two irreducible components
  $M^{\hf_1}_0$ and $M^{\hf_1}_1$, where $\dim M^{\hf_1}_0 = 24$,
  $\dim M^{\hf_1}_1 = 22$, $I_1 \in M^{\hf_1}_1 \setminus
  M^{\hf_1}_0$, and $I_0' \in M^{\hf_1}_0 \cap M^{\hf_1}_1$.
\end{theorem}

\begin{proof}
  Suppose that $J \in M^{\hf_1}$ is a point with $\gin J = I_1$.  As
  $I_1$ is saturated, so is $J$ \cite[Theorem~2.30]{Green--2010}, and
  it follows that there is a bijective morphism $J \mapsto \Proj S/J$
  between the family of such ideals and the regularity $5$ locus $R_5
  \subset H$; the locus $R_5$ is $22$-dimensional \cite[\S
    2.3]{Gotzmann--2008} and we denote by $M^{\hf_1}_1$ the closure of
  the corresponding family in $M^{\hf_1}$.  Generically, the points of
  any irreducible component of $M^{\hf_1}_1$ specialize to $I_1$.  On
  the other hand, repeating the steps in Example~\ref{defeg} shows
  that the obstruction ideal for $I_1$ is prime, meaning that
  $M^{\hf_1}$ is, analytically locally at $I_1$, an irreducible
  $22$-dimensional quadratic singularity with $\dim_{\kk} T_{I_1}
  M^{\hf_1} = 24$ \cite{Ilten--2012}.  That is, $I_1$ lies on a unique
  irreducible component of $M^{\hf_1}$ and so $M^{\hf_1}_1$ must be
  this component.  Note that the ideal $K_1 := \llrr{ x^2 - xt, xy,
    xz, y^4z, y^5 } \in M^{\hf_1}_1$ is nonsingular and specializes to
  $I_1$.

  Now suppose $J \in M^{\hf_1}$ has $\gin J = I_0'$, so that $J$ has
  regularity $6$ and is non-saturated.  Then $J^{sat}$ is in
  $M^{\hf_0}$, so that $J^{sat} = \llrr{\ell, f, g}$ is as described
  in Proposition~\ref{h0}.  Comparing $J$ with $J^{sat}$ shows that $J
  = \llrr{ \ell\ell_1, \ell\ell_2, \ell\ell_3, \ell\ell_4^4, f, g }$,
  where $\ell_1, \ell_2, \ell_3, \ell_4$ is a basis for $S_1$.  That
  is, the map $J \mapsto J^{sat}$ is a bundle over $M^{\hf_0}$ with
  fibres isomorphic to $G(3, S_1) \cong \PP^1$.  So these ideals form
  an irreducible, nonsingular, projective family of dimension $24 = 21
  + 3$, which we denote by $M^{\hf_1}_0$.  The regularity $6$ ideal
  $K_0 := \llrr{ x^2 - xt, xy, xz, xt^4, y^5, y^4z^2 } \in
  M^{\hf_1}_0$ has $\dim_{\kk} T_{K_0} M^{\hf_1} = 24$, showing that
  $M^{\hf_1}_0$ is an irreducible component of $M^{\hf_1}$.  Moreover
  the regularity $5$ ideal $\llrr{ x^2, xy, xz, xt^4 - y^4z, y^5 } =
  \llrr{ x^2, xy, xz, xt^4 - y^4z, y^5, y^4z^2}$ specializes to
  $I_0'$, showing that $I_0'$ lies in $M^{\hf_1}_0 \cap M^{\hf_1}_1$.
  In fact, the universal deformation space at $I_0'$ can be computed
  in \emph{Macaulay2}, explicitly showing that $M^{\hf_1}$ is
  analytically locally a transverse intersection of two rational
  components, of dimensions $24$ and $22$, at the point $I_0'$
  \cite{Ilten--2012}.

  Finally, if $J \in M^{\hf_1}$ is an arbitrary point, then $\gin J$
  must be either $I_0'$ or $I_1$, so $J$ belongs to either
  $M^{\hf_1}_0$ or $M^{\hf_1}_1$.
\end{proof}

\subsection{A Lex-Segment Ideal on Three Components}

The last standard-graded Hilbert scheme in this sequence is
$M^{\hf_2}$, where $\hf_2 = \hf_3 = (0, 0, 2, 8, 19, 36, 60, \ldots)$
is the Hilbert function of the ideals $I_2$ and $I_3$.  Again, we
begin with a list of its Borel-fixed ideals:
\begin{align*}
  I_0' &= \llrr{ x^2, xy, xz^2, xzt^2, xt^4, y^5, y^4z^2 },  \\
  I_1' &= \llrr{ x^2, xy, xz^2, xzt^2, y^5, y^4z },  \\
  I_2 &= \llrr{ x^2, xy, xz^2, y^4 },  \\
  I_3 &= \llrr{ x^2, xy, y^3 },
\end{align*}
derived from the saturated Borel-fixed ideals with Hilbert polynomial
$4d$ using Lemma~\ref{nonsat}.

\subsubsection{In a neighbourhood of $I_3$}

Let $J \in M^{\hf_2}$ have $\gin J = I_3$, so that $J$ is saturated,
of regularity $3$.  Then $\Proj S/J$ lies in the regularity $3$ locus
$R_3 \subset H$ \cite[\S 2.1]{Gotzmann--2008} and $J$ lies in the
corresponding $16$-dimensional family of ideals, whose closure we
denote $M^{\hf_2}_3$.  One checks that $\dim_{\kk} T_{I_3}M^{\hf_2} =
\dim_{\kk} \Hom_S(I_3, S/I_3)_0 = 16$, so $M^{\hf_2}_3$ is an
irreducible component of $M^{\hf_2}$.  Alternatively, one may observe
that $I_3$ defines an arithmetically Cohen--Macaulay scheme \cite[\S
  8]{Hartshorne--2010} that is nonsingular on $\HG$
\cite{Ellingsrud--1975} and the Comparison Theorem \cite[\S
  3]{Piene--Schlessinger--1985} holds at $I_3$.  This component
corresponds to Gotzmann's $H_{AV}$.

\subsubsection{In a neighbourhood of $I_2$}

Let $J \in M^{\hf_2}$ have $\gin J = I_2$, so that $J$ is saturated,
of regularity $4$.  So $\Proj S/J$ is in the locus $R_4 \subset H$ and
$J \in M^{\hf_2}_2$, where $M^{\hf_2}_2$ denotes the closure of the
$23$-dimensional family of ideals corresponding to $R_4$.  As $R_4$ is
irreducible, so is $M^{\hf_2}_2$.  The regularity $3$ ideal $\llrr{
  x^2, xy, xz^2-y^3 } = \llrr{ x^2, xy, xz^2-y^3, y^4 }$ specializes
to $I_2$, revealing that $I_2 \in M^{\hf_2}_2 \cap M^{\hf_2}_3$.
Notice that $K_2 := \llrr{ x^2-xt, xy, xz^2, y^4 }$ belongs to
$M^{\hf_2}_2$ and has $\dim_{\kk} T_{K_2} M^{\hf_2} = 23$, so is a
nonsingular point.  One can verify that $I_2$ satisfies the Comparison
Theorem and its universal deformation space is a transverse
intersection of a rational $16$-dimensional and a rational
$23$-dimensional component \cite{Ilten--2012}.

\subsubsection{In a neighbourhood of $I_1'$}

Let $J \in M^{\hf_2}$ have $\gin J = I_1'$, so that $J$ is
non-saturated, of regularity $5$.  As $\Proj S/J \in R_5 \subset H$,
we have $J^{sat} \in M^{\hf_1}_1$ from Theorem \ref{2comps}.  Gotzmann
shows that points of $R_5$ have saturated ideals of the form $\ell
\llrr{\ell_1, \ell_2, \ell_3} + \llrr{f_1,f_2}$, where $\ell, \ell_1,
\ell_2, \ell_3$ are linear, $\ell_1, \ell_2, \ell_3$ are linearly
independent, and $f_1,f_2 \in S_5$ (and satisfy further conditions).
Let $\ell_4$ complete the basis of $S_1$.  Then $J$ is obtained
degree-by-degree from $J^{sat}$ by: first, choosing a plane in
$J^{sat}_2 = \kk\llrr{\ell\ell_1, \ell\ell_2, \ell\ell_3}$ (say
$\kk\llrr{\ell\ell_1, \ell\ell_2}$, by a coordinate-change); second,
choosing a line in $\kk\llrr{\ell\ell_3^2, \ell\ell_3\ell_4}$; and
third, setting $J_d = J^{sat}_d$, for $d \ge 4$.  This means the map
$J \mapsto J^{sat}$ is a $\PP^1$-bundle over a $\PP^2$-bundle over the
irreducible family $M^{\hf_1}_1$, i.e.\ these ideals form an
irreducible family of dimension $22 + 2 + 1 = 25$, whose closure in
$M^{\hf_2}$ we denote by $M^{\hf_2}_1$.  A nonsingular member of
$M^{\hf_2}_1$ specializing to $I_1'$ is $K_1 := \llrr{ x^2-xt, xy,
  xz^2, xzt^2, y^5, y^4z }$.  The regularity $4$ ideal $\llrr{ x^2,
  xy, xz^2, xzt^2 - y^4 } = \llrr{ x^2, xy, xz^2, xzt^2 - y^4, y^5,
  y^4z }$ specializes to $I_1'$, showing that $I_1' \in M^{\hf_2}_2$.
Hence, we have $I_1' \in M^{\hf_2}_1 \cap M^{\hf_2}_2$.  In fact, the
power series ansatz also explicitly shows that analytically locally
$I_1'$ lies in the intersection of a rational $23$-dimensional
component and a rational $25$-dimensional component
\cite{Ilten--2012}.

\subsubsection{In a neighbourhood of $I_0'$}

Let $J \in M^{\hf_2}$ have $\gin J = I_0'$, so that $J$ is
non-saturated, of regularity $6$.  As $\Proj S/J \in R_6 \subset H$,
we know that $J^{sat} = \llrr{\ell, \ell' h, qh} \in M^{\hf_0}$ has
the form described in Proposition~\ref{h0}.  Comparing Hilbert
functions shows that $J$ is obtained from $J^{sat}$ as follows: first,
remove $\ell$; second, choose a plane in $J^{sat}_2 = S_1\ell$ (say
$\kk\llrr{\ell x, \ell y}$ up to a change-of-basis); third, as
$\dim_{\kk} S_1 J_2 = 7$, choose a line in $J^{sat}_3 \setminus S_1
J_2 = S_2\ell \setminus S_1 J_2$ (say spanned by some $\ell q' \in
\kk\llrr{\ell z^2, \ell zt, \ell t^2}$); fourth, as $\dim_{\kk} S_1
J_3 = 18$, choose a line in $J^{sat}_4 \setminus S_1 J_3 = S_3\ell
\setminus S_1 J_3$ (say spanned by some $\ell p \in \kk\llrr{ \ell
  z^3, \ell z^2t, \ell zt^2, \ell t^3 } \setminus \kk\llrr{\ell q' z,
  \ell q' t}$); finally, set $J_d = J^{sat}_d$, for $d \ge 5$.  In
other words, the map $J \mapsto J^{sat}$ is a $\PP^1$-bundle over a
$\PP^2$-bundle over a $G(2,S_1\ell)$-bundle over $M^{\hf_0}$.  This
shows that the family of such ideals is nonsingular, irreducible,
projective, and of dimension $21 + 4 + 2 + 1 = 28$.  Of course, the
lex-segment ideal $I_0'$ is one such ideal.  Another such ideal is
$K_0 := \llrr{x^2, xy-xt, xz^2, xzt^2, xt^4, y^5, y^4z^2 }$, which has
$\dim_{\kk} T_{K_0} M^{\hf_2} = 28$, so the family is an irreducible
component, which we denote by $M^{\hf_2}_0$.

Consider the nonsingular point $K_2 := \llrr{ x^2-xt, xy, xz^2, y^4 }
\in M^{\hf_2}_2$.  The change-of-basis $\gamma \colon x \mapsto x+t, y
\mapsto x+y, z \mapsto x, t \mapsto z$ transforms $K_2$ into the ideal
\begin{align*}
  \gamma K_2 := \langle & xy + xz - xt + yt + zt - t^2, \quad x^2 - xz
  + 2xt - zt + t^2, \\ & \qquad xz^2 - 2xzt + z^2t + xt^2 - 2zt^2 +
  t^3, \\ & \qquad \qquad y^4 - 4y^3t + 6y^2t^2 - 3xzt^2 + 2xt^3 -
  4yt^3 - 3zt^3 + 3t^4 \rangle
\end{align*}
in $M^{\hf_2}_2$ and the lexicographic initial ideal of $\gamma K_2$
is $I_0'$, showing that $I_0'$ also belongs to $M^{\hf_2}_2$.
Moreover, the ideal $\llrr{x^2, xy, xz^2, xzt^2, xt^4 - y^4z, y^5} =
\llrr{x^2, xy, xz^2, xzt^2, xt^4 - y^4z, y^5, y^4z^2}$ has regularity
$5$ and specializes to $I_0'$, showing that $I_0' \in M^{\hf_2}_1$.
Hence, we have shown that $I_0'$ lies in the triple intersection
$M^{\hf_2}_0 \cap M^{\hf_2}_1 \cap M^{\hf_2}_2$.  To see that $I_0'
\notin M^{\hf_2}_3$, note that any subfamily of $M^{\hf_2}_3$
containing $I_0'$ defines a flat family in $H_{AV}$ containing the
lexicographic point $\Proj S/I_0'$; this is a contradiction by
smoothness of $\Proj S/I_0'$ \cite{Reeves--Stillman--1997}.  The power
series ansatz terminates for $I_0'$, but the obstruction ideal is too
complicated to directly obtain its primary decomposition.  However,
one quickly identifies a unique $28$-dimensional irreducible component
and the residual scheme has dimension $25$, as expected from our
description.

\begin{theorem}
  The standard-graded Hilbert scheme $M^{\hf_2}$ consists of four
  irreducible components, of dimensions $16, 23, 25$, and $28$.  The
  lex-segment ideal with Hilbert function $\hf_2$ lies in the
  intersection of the three highest-dimensional components.
\end{theorem}

\begin{proof}
  If $J \in M^{\hf_2}$, then $\gin J$ is one of the aforementioned
  Borel-fixed ideals, so $J$ belongs to one of the identified
  components $M^{\hf_2}_0, M^{\hf_2}_1, M^{\hf_2}_2$, or
  $M^{\hf_2}_3$.  Hence, the preceding paragraphs identify all
  irreducible components of $M^{\hf_2}$.  (In particular, a
  nonsingular point is identified in each component.)
\end{proof}

\begin{remark}
  The incidence complex of $M = M^{\hf_2} =
  \hilb^{\hf_2}(\kk[x,y,z,t])$ is as follows.  Here the upper index at
  a node displays the dimension of the corresponding irreducible
  component.
  
  
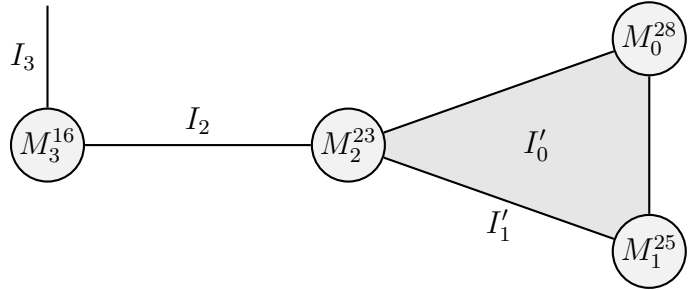
{%
  \centering
  \begin{tikzpicture}
    \tikzstyle{point}=[circle,thick,draw=black,fill=gray!10,inner
      sep=1pt,minimum width=4pt,minimum height=4pt]
    
    \node (a)[point] at (0,0) {$M^{16}_3$};
    \node (b)[point] at (4,0) {$M^{23}_2$};
    \node (c)[point] at (8,-2/1.414) {$M^{25}_1$};
    \node (d)[point] at (8,2/1.414) {$M^{28}_0$};    
    \node (e)[] at (6.5,0) {$I_0'$};
    \node (o)[] at (0,2) {};

    \draw[thick] (a) -- node[left] {$I_3$} (o);
    \draw[thick] (a) -- node[above] {$I_2$} (b);
    \draw[thick] (b) -- node[below] {$I_1'$} (c) -- (d) --
    (b);

    \begin{pgfonlayer}{background}
      \fill[color=gray!20] (b.center) -- (c.center) -- (d.center) --
      cycle;
    \end{pgfonlayer}
  \end{tikzpicture}
  
  }  
\end{remark}
% --------------------------------------------------------------------

% --------------------------------------------------------------------
% Section 3 -- Further Examples
% --------------------------------------------------------------------
\section{Nearby Examples}

The analysis of the previous section can be performed on many
examples.  Here we do this at the nodes immediately preceding
Gotzmann's example in the codimension $2$ Hilbert tree.

\subsection{An Example in the Plane}

The Hilbert scheme $\hilb^4(\PP^2)$ has two Borel-fixed points, with
saturated ideals in $S := \kk[x,y,z]$ and Hilbert functions as
follows:
\begin{align*}
  I_0 &= \llrr{ x, y^4 }, &\hf_0 &= (0, 1, 3, 6, 11, 17, \ldots), \\
  I_1 &= \llrr{ x^2, xy, y^3 }, &\hf_1 &= (0, 0, 2, 6, 11, 17,
  \ldots),
\end{align*}
where $I_i$ has regularity $4 - i$.  The standard-graded Hilbert
scheme $M^{\hf_0} := \hilb^{\hf_0}(S)$ is nonsingular and irreducible,
of dimension $6$, comprised of ideals of the form $\llrr{\ell, f}$
where $\ell \in S_1$, $f \in S_4$, and $\ell \nmid f$.

We study the standard-graded Hilbert scheme $M^{\hf_1} :=
\hilb^{\hf_1}(S)$.  Its Borel-fixed points are
\begin{align*}
  I_0' &= \llrr{ x^2, xy, xz^2, y^4 }, \\
  I_1 &= \llrr{ x^2, xy, y^3 },
\end{align*}
where $I_0'$ is obtained from $I_0$ via Lemma~\ref{nonsat}.  We prove
the following result.

\begin{theorem}
  The standard-graded Hilbert scheme $M^{\hf_1}$ is a union of two
  $8$-dimensional components $M^{\hf_1}_0, M^{\hf_1}_1$ with
  lex-segment ideal $I_0' \in M^{\hf_1}_0 \cap M^{\hf_1}_1$ and $I_1
  \in M^{\hf_1}_1 \setminus M^{\hf_1}_0$.
\end{theorem}

\begin{proof}
  Every $J \in M^{\hf_1}$ has $J_2 = \kk\llrr{q_1, q_2}$ spanned by
  two quadrics.  If $q_1, q_2$ do not share a common factor, then $J =
  \llrr{q_1, q_2}$ is a complete intersection ideal and the Koszul
  complex shows that $J$ has regularity $3$; moreover, $J$ is
  saturated and $\gin J = I_1$.  Conversely, suppose that $J$ has
  $\gin J = I_1$, so that $J$ is saturated of regularity $3$, and
  suppose that $q_1, q_2$ do share a common factor, say $q_1 =
  \ell\ell_1, q_2 = \ell\ell_2$, where $\ell, \ell_1, \ell_2 \in S_1$
  and $\ell_1, \ell_2$ are linearly independent.  Then there is a
  cubic generator $p \in J_3 \setminus S_1J_2$ and we must have $\ell
  \nmid p$, as $J$ is saturated but $\ell \notin J$.  This shows that
  $J = \llrr{\ell\ell_1, \ell\ell_2, p}$ with $\ell \nmid p$;
  moreover, $\hf_1(4) = 11$ implies $p \in \llrr{\ell_1, \ell_2}$.
  Let $M^{\hf_1}_1$ denote the closure of the family of ideals $J$
  with $\gin J = I_1$.  Then the subfamily of complete intersections
  has dimension $8$, while the family of ideals with a cubic generator
  has dimension $7$, so $M^{\hf_1}_1$ has dimension $8$.  In fact, we
  have $\dim_{\kk} T_{I_1} M^{\hf_1} = 8$, so $I_1$ is a nonsingular
  point and $M^{\hf_1}_1$ is an irreducible component of $M^{\hf_1}$.

  On the other hand, if $J$ has $\gin J = I_0'$, then $J$ is
  non-saturated, of regularity $4$, and $J^{sat} \in M^{\hf_0}$ has
  the form $\llrr{\ell, f}$ where $\ell \in S_1$, $f \in S_4$, and
  $\ell \nmid f$.  To obtain $J$ from $J^{sat}$, first, remove $\ell$,
  second, choose a plane $\kk\llrr{\ell\ell_1, \ell\ell_2} \subset
  J^{sat}_2 = S_1\ell$, and third, set $J_d = J^{sat}_d$, for $d \ge
  3$.  This family is a $\PP^2$-bundle over $M^{\hf_0}$, giving an
  irreducible, nonsingular, projective family of dimension $6 + 2 =
  8$, which we denote by $M^{\hf_1}_0$.  The ideal $K_0 :=
  \llrr{x^2-xz, xy, xz^2, y^4}$ is in $M^{\hf_1}_0$, with $\dim_{\kk}
  T_{K_0} M^{\hf_1} = 8$, so $M^{\hf_1}_0$ is an irreducible
  component.  Moreover, the ideal $\llrr{x^2, xy, xz^2 - y^3} =
  \llrr{x^2, xy, xz^2 - y^3, y^4}$ has regularity $3$ (and is
  nonsingular on $M^{\hf_1}_1$) but specializes to $I_0'$.  This shows
  that $I_0' \in M^{\hf_1}_0 \cap M^{\hf_1}_1$.

  Any $J \in M^{\hf_1}$ has $\gin J = I_0'$ or $\gin J = I_1$, so
  $M^{\hf_1} = M^{\hf_1}_0 \cup M^{\hf_1}_1$.  In fact, the power
  series ansatz at $I_0'$ yields a transverse intersection of two
  $8$-dimensional components.
\end{proof}

\subsection{Lifting to Degree Four Space Curves}

The Hilbert scheme $\hilb^{4d-2}(\PP^3)$ has a unique Borel-fixed
point, so is nonsingular and irreducible \cite[Lemma
  5.6]{Staal--2020}.  If $\hf_0$ denotes its Hilbert function and
$I_0$ its ideal, then every homogeneous ideal $J \subset S =
\kk[x,y,z,t]$ with Hilbert function $\hf_0$ has $\gin J = I_0$ and is
saturated.  By Proposition~\ref{h0'}, the standard-graded Hilbert
scheme $\hilb^{\hf_0}(S)$ is nonsingular and irreducible; its
dimension is $17$ and its points are complete intersection ideals
$\llrr{\ell, f}$, for $\ell \in S_1$ and $f \in S_4$.

Moving towards Gotzmann's example, the Hilbert scheme
$\hilb^{4d-1}(\PP^3)$ has two Borel-fixed points, with saturated
ideals in $S = \kk[x,y,z,t]$ and Hilbert functions as follows:
\begin{align*}
  I_0 &= \llrr{ x,y^5,y^4z }, &\hf_0 &= (0, 1, 4, 10, 20, 37, 61,
  \ldots), \\
  I_1 &= \llrr{ x^2, xy, xz, y^4 }, &\hf_1 &= (0, 0, 3, 9, 20, 37,
  61, \ldots),
\end{align*}
where $I_i$ has regularity $5-i$.  As in Proposition~\ref{h0}, the
standard-graded Hilbert scheme $\hilb^{\hf_0}(S)$ is nonsingular and
irreducible; its dimension is $19$ and its points have the form
$\llrr{\ell, \ell_1 g, \ell_2 g}$, where $\ell, \ell_1, \ell_2 \in
S_1$ are linearly independent and $g \in S_4$ such that $\ell \nmid
g$.

So let us examine $M^{\hf_1} := \hilb^{\hf_1}(S)$.  Its Borel-fixed
ideals are
\begin{align*}
  I_0' &= \llrr{ x^2, xy, xz, xt^3, y^5, y^4z }, \\
  I_1 &= \llrr{ x^2, xy, xz, y^4 },
\end{align*}
where $I_0'$ is obtained from $I_0$ via Lemma~\ref{nonsat}.  We aim to
prove the following.

\begin{theorem}
  \label{4d-1}
  We have $M^{\hf_1} = M^{\hf_1}_0 \cup M^{\hf_1}_1$, where
  $M^{\hf_1}_0$ is a $22$-dimensional component, $M^{\hf_1}_1$ is a
  $20$-dimensional component, the lex-segment ideal $I_0' \in
  M^{\hf_1}_0 \cap M^{\hf_1}_1$, and $I_1 \in M^{\hf_1}_1 \setminus
  M^{\hf_1}_0$.
\end{theorem}

\subsubsection{In a neighbourhood of $I_1$}

Suppose that $J \in M^{\hf_1}$ has $\gin J = I_1$, so that $J$ is
saturated, of regularity $4$.  We have $J_2 = \kk\llrr{q_1, q_2,
  q_3}$, where the quadrics do not generate a complete intersection
ideal, so at least two of the quadrics share a common factor, say $q_1
= \ell \ell_1$ and $q_2 = \ell \ell_2$ with $\ell_1, \ell_2$ linearly
independent.  But then $\dim_{\kk} S_1 \llrr{q_1, q_2} \cap S_1 q_3 =
2$, which is only possible if $\ell \vert q_3$ as well, so we set $q_3
= \ell \ell_3$.  We find that $\dim_{\kk} S_1 J_2 = 9$ and $\dim_{\kk}
S_2 J_2 = 19$, so that $J_3 = S_1 J_2$ and there exists $g \in J_4
\setminus S_2 J_2$.  If $\ell \vert g$, then without loss of
generality $g = \ell \ell_4^3$, where $\ell_1, \ell_2, \ell_3, \ell_4$
is a basis of $S_1$, but then $J = J^{sat}$ implies $\ell \in J$, a
contradiction.  Thus, we have $J = \llrr{\ell\ell_1, \ell\ell_2,
  \ell\ell_3, g}$ with $g \in S_4$, $\ell \nmid g$, and moreover $g
\in \llrr{\ell_1, \ell_2, \ell_3}$.  A parameter count shows that this
family has dimension $20$; we denote its Zariski closure by
$M^{\hf_1}_1$.  One verifies that $\dim_{\kk} T_{I_1} M^{\hf_1} = 20$,
so $M^{\hf_1}_1$ is an irreducible component with nonsingular point
$I_1$.

\subsubsection{In a neighbourhood of $I_0'$}

If $J \in M^{\hf_1}$ has $\gin J = I_0'$, then $J$ is non-saturated,
of regularity $5$.  This means $J^{sat} \in M^{\hf_0}$, so $J^{sat} =
\llrr{\ell, \ell_1' g, \ell_2' g}$, where $\ell, \ell_1', \ell_2' \in
S_1$ are linearly independent and $g \in S_4$ such that $\ell \nmid
g$.  We obtain $J$ from $J^{sat}$ degree-by-degree by, first, removing
$\ell$, and second, choosing a hyperplane $\kk\llrr{\ell\ell_1,
  \ell\ell_2, \ell\ell_3}$ in $J^{sat}_2 = S_1\ell$; we then have $J_3
= S_1 J_2$ and $J_d = J^{sat}_d$, for $d \ge 4$.  This means the
family of such ideals has the structure of a $\PP^3$-bundle over
$M^{\hf_0}$, resulting in a nonsingular, irreducible, projective
family of dimension $19 + 3 = 22$, which we denote by $M^{\hf_1}_0$.
The ideal $K_0 = \llrr{x^2 - xt, xy, xz, xt^3, y^5, y^4z}$ belongs to
$M^{\hf_1}_0$ and has $\dim_{\kk} T_{K_0} M = 22$, so $M^{\hf_1}_0$ is
an irreducible component.  Further, the regularity $4$ ideal
$\llrr{x^2, xy, xz, xt^3 - y^4} = \llrr{x^2, xy, xz, xt^3 - y^4, y^5,
  y^4z} \in M^{\hf_1}_1$ specializes to $I_0'$, showing that the
lex-segment ideal satisfies $I_0' \in M^{\hf_1}_0 \cap M^{\hf_1}_1$.

\begin{proof}[Proof of Theorem~\ref{4d-1}]
  If $J \in M^{\hf_1}$, then either $\gin J = I_0'$ or $\gin J = I_1$,
  so $J$ belongs to one of the identified components.  Note that the
  power series ansatz terminates at $I_0'$, showing the transverse
  intersection of $M^{\hf_1}_0$ and $M^{\hf_1}_1$ \cite{Ilten--2012}.
\end{proof}

% --------------------------------------------------------------------
% Section 5 -- Twisted Cubics
% --------------------------------------------------------------------
\section{Twisted Cubics}

We conclude with an example that has been hiding in plain sight for
some time.

\subsection{A Singular Lex-Segment Ideal of Twisted Cubics}

The Borel-fixed points on $\hilb^{3d+1}(\PP^3)$ have saturated ideals
in $S = \kk[x,y,z,t]$ with Hilbert functions as follows:
\begin{align*}
  I_0 &= \llrr{ x, y^4, y^3z }, &\hf_0 &= (0, 1, 4, 10, 22, 40, 65,
  \ldots), \\
  I_1 &= \llrr{ x^2, xy, xz, y^3 }, &\hf_1 &= (0, 0, 3, 10, 22, 40,
  65, \ldots), \\
  I_2 &= \llrr{ x^2, xy, y^2 }, &\hf_2 &= (0, 0, 3, 10, 22, 40,
  65, \ldots),
\end{align*}
where $I_i$ has regularity $4-i$.

Similarly to previous instances, $M^{\hf_0} := \hilb^{\hf_0}(S)$ is a
nonsingular and irreducible projective variety; its dimension is $14$
and its points have the form $\llrr{\ell, \ell_1 p, \ell_2 p}$, where
$\ell, \ell_1, \ell_2$ are linearly independent, $p \in S_3$, and
$\ell \nmid p$.

Next, as $\hf_1 = \hf_2$, we study $M^{\hf_1} := \hilb^{\hf_1}(S)$.
Its Borel-fixed points are the following:
\begin{align*}
  I_0' &= \llrr{ x^2, xy, xz, xt^2, y^4, y^3z }, \\
  I_1 &= \llrr{ x^2, xy, xz, y^3 }, \\
  I_2 &= \llrr{ x^2, xy, y^2 },
\end{align*}
where $I_0'$ is obtained via Lemma~\ref{nonsat}.  The analysis is
analogous to that of previous examples.

\subsubsection{In a neighbourhood of $I_2$}

If $J \in M^{\hf_1}$ has $\gin J = I_2$, then $J$ is saturated, of
regularity $2$, and so $J = \llrr{q_1, q_2, q_3}$ is generated by
quadrics.  This family parametrizes the twisted cubics.  We denote its
closure by $M^{\hf_2}_2$ and note that $\dim M^{\hf_1}_2 = \dim_{\kk}
T_{I_2} M^{\hf_1} = 12$, so that $M^{\hf_1}_2$ is an irreducible
component with nonsingular point $I_2$.

\subsubsection{In a neighbourhood of $I_1$}

If $J \in M^{\hf_1}$ has $\gin J = I_1$, then $J$ is saturated, of
regularity $3$, and again we know that $J_2 = \kk\llrr{q_1, q_2,
  q_3}$, where $q_1, q_2, q_3$ now share a common factor.  So $q_1 =
\ell\ell_1, q_2 = \ell\ell_2, q_3 = \ell\ell_3$, for some linear
$\ell$ and linearly independent $\ell_1, \ell_2, \ell_3 \in S_1$.
These quadrics generate a $9$-dimensional subspace of $S_3$, so there
is an additional cubic generator, say $p \in J_3$.  Then $\ell \nmid
p$, because otherwise $J = J^{sat}$ implies $\ell \in J$; we further
see that $\dim _{\kk} S_1 p \cap S_2 \llrr{\ell \ell_1, \ell \ell_2,
  \ell \ell_3} = 1$, so $p \in \llrr{\ell_1, \ell_2, \ell_3}$.  This
family is in bijection with a dense open subset of the lexicographic
component of $\hilb^{3d+1}(\PP^3)$, in other words, is an irreducible
component of dimension $15$ \cite[Lemma~4]{Piene--Schlessinger--1985};
we denote its closure by $M^{\hf_1}_1$.  The ideal $\llrr{x^2, xy, xz
  - y^2} = \llrr{x^2, xy, xz - y^2, y^3}$ specializes to $I_1$ but has
regularity $2$, showing that $I_1 \in M^{\hf_1}_1 \cap M^{\hf_1}_2$.
This is expected of course, as $I_1$ satisfies the Comparison Theorem
and the well-known description of the component structure of
$\hilb^{3d+1}(\PP^3)$ hinges on the termination of the power series
ansatz at $I_1$ \cite[Lemma 6]{Piene--Schlessinger--1985}.

\subsubsection{In a neighbourhood of $I_0'$}

If $J \in M^{\hf_1}$ has $\gin J = I_0'$, then $J$ is non-saturated,
of regularity $4$, and $J^{sat} \in M^{\hf_0}$ equals $\llrr{\ell,
  \ell_1' p, \ell_2' p}$, where $\ell, \ell_1', \ell_2'$ are linearly
independent, $p \in S_3$, and $\ell \nmid p$.  We obtain $J$ from
$J^{sat}$ by first, removing $\ell$, second, choosing a hyperplane
$\kk\llrr{\ell\ell_1, \ell\ell_2, \ell\ell_3} \subset J^{sat}_2 =
S_1\ell$, and third, setting $J_d = J^{sat}_d$, for $d \ge 3$.  The
family of such ideals forms a $\PP^3$-bundle over $M^{\hf_0}$,
defining an irreducible, nonsingular, projective family of dimension
$14 + 3 = 17$, which we denote by $M^{\hf_1}_0$.  The ideal $K_0 =
\llrr{x^2 - xt, xy, xz, xt^2, y^4, y^3z}$ is in $M^{\hf_1}_0$ and has
$\dim_{\kk} T_{K_0} M^{\hf_1} = 17$, showing that $M^{\hf_1}_0$ is an
irreducible component of $M^{\hf_1}$.  Observe that the ideal
$\llrr{x^2, xy, xz, xt^2 - y^3} = \llrr{x^2, xy, xz, xt^2 - y^3, y^3z,
  y^4}$ has regularity $3$ and specializes to $I_0'$, showing that
$I_0' \in M^{\hf_1}_1$.  Hence, the lex-segment ideal $I_0'$ lies in
the intersection $M^{\hf_1}_0 \cap M^{\hf_1}_1$.  In fact,
Example~\ref{defeg} demonstrates that the power series ansatz
terminates at $I_0'$, showing that the universal deformation space is
a transverse intersection of a rational $15$-dimensional component and
a rational $17$-dimensional component \cite{Ilten--2012}.

\begin{theorem}
  \label{twist}
  We have $M^{\hf_1} = M^{\hf_1}_0 \cup M^{\hf_1}_1 \cup M^{\hf_1}_2$,
  where the irreducible components have dimensions $\dim M^{\hf_1}_0 =
  17$, $\dim M^{\hf_1}_1 = 15$, and $\dim M^{\hf_1}_2 = 12$, and the
  lex-segment ideal lies in $M^{\hf_1}_0 \cap M^{\hf_1}_1$.
\end{theorem}

\begin{proof}
  Any ideal $J$ with Hilbert function $\hf_1$ must have $\gin J$ equal
  to $I_0'$, $I_1$, or $I_2$, and hence lies in one of the identified
  components.  The preceding paragraphs derive the decomposition and
  show $I_0' \in M^{\hf_1}_0 \cap M^{\hf_1}_1$.
\end{proof}

\begin{remark}
  The incidence graph of $M = M^{\hf_1}$ is as follows.  The upper
  index at a node displays the dimension of the corresponding
  irreducible component.

  
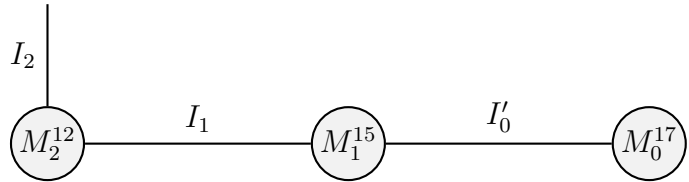
{%
  \centering
  \begin{tikzpicture}
    \tikzstyle{point}=[circle,thick,draw=black,fill=gray!10,inner sep=1pt,minimum
      width=4pt,minimum height=4pt]
    
    \node (a)[point] at (0,0) {$M^{12}_2$};
    \node (b)[point] at (4,0) {$M^{15}_1$};
    \node (c)[point] at (8,0) {$M^{17}_0$};
    \node (o)[] at (0,2) {};

    \draw[thick] (a) -- node[left] {$I_2$} (o);
    \draw[thick] (a) -- node[above] {$I_1$} (b);
    \draw[thick] (b) -- node[above] {$I_0'$} (c);
  \end{tikzpicture}
  
  }  
\end{remark}
% --------------------------------------------------------------------

% --------------------------------------------------------------------
% Bibliography
% --------------------------------------------------------------------
 \bibliography{curves}{}

\providecommand{\bysame}{\leavevmode\hbox to3em{\hrulefill}\thinspace}
\providecommand{\MR}{\relax\ifhmode\unskip\space\fi MR }
\providecommand{\MRhref}[2]{%
  \href{http://www.ams.org/mathscinet-getitem?mr=#1}{#2}
}
\providecommand{\href}[2]{#2}
\begin{thebibliography}{{Got}08}

\bibitem[Bay82]{Bayer--1982}
David~Allen Bayer, \emph{The division algorithm and the {H}ilbert scheme},
  ProQuest LLC, Ann Arbor, MI, 1982, Thesis (Ph.D.)--Harvard University.
  \MR{2632095}

\bibitem[BH93]{Bruns--Herzog--1993}
Winfried Bruns and J{\"u}rgen Herzog, \emph{Cohen-{M}acaulay rings}, Cambridge
  Studies in Advanced Mathematics, vol.~39, Cambridge University Press,
  Cambridge, 1993. \MR{1251956 (95h:13020)}

\bibitem[BS87a]{Bayer--Stillman--1987}
David Bayer and Michael Stillman, \emph{A criterion for detecting
  {$m$}-regularity}, Invent. Math. \textbf{87} (1987), no.~1, 1--11. \MR{862710
  (87k:13019)}

\bibitem[BS87b]{Bayer--Stillman--1987b}
\bysame, \emph{A theorem on refining division orders by the reverse
  lexicographic order}, Duke Math. J. \textbf{55} (1987), no.~2, 321--328.
  \MR{894583}

\bibitem[CLO15]{Cox--Little--O'Shea--2015}
David~A. Cox, John Little, and Donal O'Shea, \emph{Ideals, varieties, and
  algorithms}, fourth ed., Undergraduate Texts in Mathematics, Springer, Cham,
  2015, An introduction to computational algebraic geometry and commutative
  algebra. \MR{3330490}

\bibitem[Eis95]{Eisenbud--1995}
David Eisenbud, \emph{Commutative algebra}, Graduate Texts in Mathematics, vol.
  150, Springer-Verlag, New York, 1995, With a view toward algebraic geometry.
  \MR{1322960}

\bibitem[Ell75]{Ellingsrud--1975}
Geir Ellingsrud, \emph{Sur le sch\'{e}ma de {H}ilbert des vari\'{e}t\'{e}s de
  codimension {$2$} dans {${\bf P}\sp{e}$} \`a c\^{o}ne de {C}ohen-{M}acaulay},
  Ann. Sci. \'{E}cole Norm. Sup. (4) \textbf{8} (1975), no.~4, 423--431.
  \MR{393020}

\bibitem[Gal74]{Galligo--1974}
Andr{\'e} Galligo, \emph{\`{A} propos du th\'eor\`eme de-pr\'eparation de
  {W}eierstrass}, Fonctions de plusieurs variables complexes ({S}\'em. {F}ran\c
  cois {N}orguet, octobre 1970--d\'ecembre 1973; \`a la m\'emoire d'{A}ndr\'e
  {M}artineau), Springer, Berlin, 1974, Th{\`e}se de 3{\`e}me cycle soutenue le
  16 mai 1973 {\`a} l'Institut de Math{\'e}matique et Sciences Physiques de
  l'Universit{\'e} de Nice, pp.~543--579. Lecture Notes in Math., Vol. 409.
  \MR{0402102}

\bibitem[{Got}08]{Gotzmann--2008}
Gerd {Gotzmann}, \emph{{The irreducible components of
  $Hilb^{4n}(\mathbb{P}^3)$}}, arXiv e-prints (2008), arXiv:0811.3160.

\bibitem[Gre10]{Green--2010}
Mark~L. Green, \emph{Generic initial ideals [mr1648665]}, Six lectures on
  commutative algebra, Mod. Birkh\"auser Class., Birkh\"auser Verlag, Basel,
  2010, pp.~119--186. \MR{2641237}

\bibitem[GS]{Grayson--Stillman}
Daniel~R. Grayson and Michael~E. Stillman, \emph{Macaulay2, a software system
  for research in algebraic geometry}, Available at
  http://www.math.uiuc.edu/Macaulay2/.

\bibitem[GW10]{Goertz--Wedhorn--2010}
Ulrich G\"{o}rtz and Torsten Wedhorn, \emph{Algebraic geometry {I}}, Advanced
  Lectures in Mathematics, Vieweg + Teubner, Wiesbaden, 2010, Schemes with
  examples and exercises. \MR{2675155}

\bibitem[Har10]{Hartshorne--2010}
Robin Hartshorne, \emph{Deformation theory}, Graduate Texts in Mathematics,
  vol. 257, Springer, New York, 2010. \MR{2583634}

\bibitem[HS04]{Haiman--Sturmfels--2004}
Mark Haiman and Bernd Sturmfels, \emph{Multigraded {H}ilbert schemes}, J.
  Algebraic Geom. \textbf{13} (2004), no.~4, 725--769. \MR{2073194}

\bibitem[Ilt12]{Ilten--2012}
Nathan~Owen Ilten, \emph{Versal deformations and local {H}ilbert schemes}, J.
  Softw. Algebra Geom. \textbf{4} (2012), 12--16. \MR{2947667}

\bibitem[Mac27]{Macaulay--1927}
F.~S. Macaulay, \emph{Some {P}roperties of {E}numeration in the {T}heory of
  {M}odular {S}ystems}, Proc. London Math. Soc. \textbf{S2-26} (1927), no.~1,
  531. \MR{1576950}

\bibitem[MN14]{Moore--Nagel--2014}
Dennis Moore and Uwe Nagel, \emph{Algorithms for strongly stable ideals}, Math.
  Comp. \textbf{83} (2014), no.~289, 2527--2552. \MR{3223345}

\bibitem[MS05]{Miller--Sturmfels--2005}
Ezra Miller and Bernd Sturmfels, \emph{Combinatorial commutative algebra},
  Graduate Texts in Mathematics, vol. 227, Springer-Verlag, New York, 2005.
  \MR{2110098 (2006d:13001)}

\bibitem[MS10]{Maclagan--Smith--2010}
Diane Maclagan and Gregory~G. Smith, \emph{Smooth and irreducible multigraded
  {H}ilbert schemes}, Adv. Math. \textbf{223} (2010), no.~5, 1608--1631.
  \MR{2592504}

\bibitem[PS85]{Piene--Schlessinger--1985}
Ragni Piene and Michael Schlessinger, \emph{On the {H}ilbert scheme
  compactification of the space of twisted cubics}, Amer. J. Math. \textbf{107}
  (1985), no.~4, 761--774. \MR{796901}

\bibitem[PS02]{Peeva--Stillman--2002}
Irena Peeva and Mike Stillman, \emph{Toric {H}ilbert schemes}, Duke Math. J.
  \textbf{111} (2002), no.~3, 419--449. \MR{1885827}

\bibitem[Ree92]{Reeves--1992}
Alyson~April Reeves, \emph{Combinatorial structure on the {H}ilbert scheme},
  ProQuest LLC, Ann Arbor, MI, 1992, Thesis (Ph.D.)--Cornell University.
  \MR{2688239}

\bibitem[RS97]{Reeves--Stillman--1997}
Alyson Reeves and Mike Stillman, \emph{Smoothness of the lexicographic point},
  J. Algebraic Geom. \textbf{6} (1997), no.~2, 235--246. \MR{1489114
  (98m:14003)}

\bibitem[RS20]{Ramkumar--Sammartano--2020}
Ritvik Ramkumar and Alessio Sammartano, \emph{Singular lexicographic points on
  hilbert schemes}, 2020.

\bibitem[San05]{Santos--2005}
Francisco Santos, \emph{Non-connected toric {H}ilbert schemes}, Math. Ann.
  \textbf{332} (2005), no.~3, 645--665. \MR{2181765}

\bibitem[SS20]{Skjelnes--Smith--2020}
Roy {Skjelnes} and Gregory~G. {Smith}, \emph{{Smooth Hilbert schemes: their
  classification and geometry}}, arXiv e-prints (2020), arXiv:2008.08938.

\bibitem[Sta]{Staal--2014}
Andrew~P. Staal, \emph{Computational {D}eformation {T}heory of {P}rojective
  {S}chemes}, Comprehensive Examination Document, Queen's University, available
  at \url{http://astaal.be/research.html}.

\bibitem[Sta20]{Staal--2020}
\bysame, \emph{The ubiquity of smooth {H}ilbert schemes}, Math. Z. \textbf{296}
  (2020), no.~3-4, 1593--1611. \MR{4159841}

\bibitem[Ste95]{Stevens--1995}
Jan Stevens, \emph{Computing versal deformations}, Experiment. Math. \textbf{4}
  (1995), no.~2, 129--144. \MR{1377414}

\end{thebibliography}
 \bibliographystyle{amsalpha}
% --------------------------------------------------------------------

\end{document}